# q−Catalan numbers and q−Narayana polynomials

Johann Cigler

Fakultät für Mathematik
Universität Wien
A-1090 Wien, Nordbergstraße 15
e-mail: johann.cigler@univie.ac.at

The q−Catalan numbers
$1, 1, 1+q, 1+2q+q^2+q^3, 1+3q+3q^2+3q^3+2q^4+q^5+q^6, \cdots,$
which were invented by Carlitz and Riordan ([2], cf. also [7] and [10]), are defined by
$C_n(q) = \sum_{k=0}^{n-1} q^k C_k(q) C_{n-k-1}(q)$ with initial value $C_0(q) = 1$.

Let $f(z,q) = \sum_{k \geq 0} C_k(q) z^k$ be their generating function, which is uniquely determined by the functional equation $f(z,q) = 1 + z f(z,q) f(qz,q)$. This implies the well known fact that it can be represented in the form

$$f(z,q) = \frac{E_2(-qz)}{E_2(-z)}, \tag{1}$$

where $E_r(z)$ denotes the generalized q−exponential function

$$E_r(z) = \sum_{k \geq 0} q^{r\binom{k}{2}} \frac{z^k}{(1-q)(1-q^2)\cdots(1-q^k)}. \tag{2}$$

In this note we want to sketch some extensions of this interesting result.

a) First observe that
$$E_r(z) - E_r(qz) = z E_r(q^r z). \tag{3}$$

For $n \in \mathbb{N}$ define

$$G_r(z,n) = \sum_{k \geq 0} G(k,n,r) z^k := \frac{E_r(-q^n z)}{E_r(-z)}. \tag{4}$$

Then we have

$$G_r(z,n+1) = G_r(z,n) + q^n z G_r(z, n+r). \tag{5}$$

Comparing coefficients we get

$$\frac{G(k,n+1,r) - G(k,n,r)}{q^n} = G(k-1, n+r, r) \tag{6}$$

with $G(k,0,r) = [k=0]$ and $G(0,n,r) = 1$.

This implies
$$G_r(z,1) = 1 + zG_r(z,r). \qquad (7)$$

These are the characteristic properties of the $q-$Gould polynomials (cf. [4]).

For $q=1$ they have the explicit formula $G(k,n,r) = \dfrac{n}{n+rk}\binom{n+rk}{k}$ (cf. e.g. [9]).

For general $q$ special values are $G(k,n,0) = q^{\binom{k}{2}}\begin{bmatrix}n\\k\end{bmatrix}$ and $G(k,n,1) = \begin{bmatrix}n+k-1\\k\end{bmatrix}$, where $\begin{bmatrix}n\\k\end{bmatrix}$
denotes a $q-$binomial coefficient. For $r>1$ no explicit formulas are known.

Note that $G_2(z,1) = 1 + zG_2(z,2) = 1 + zG_2(z,1)G_2(qz,1) = f(z,q)$ is the generating function of the $q-$Catalan numbers.

From $\dfrac{E_r(-q^n z)}{E_r(-z)} = \dfrac{E_r(-qz)}{E_r(-z)}\dfrac{E_r(-q^2 z)}{E_r(-qz)}\cdots\dfrac{E_r(-q^n z)}{E_r(-q^{n-1}z)}$

we get
$$G_r(z,n) = G_r(z,1)G_r(qz,1)\cdots G_r(q^{n-1}z,1) \qquad (8)$$

and
$$G_r(z, m+n) = G_r(z,m)G_r(q^m z, n). \qquad (9)$$

b) Let now $(a \dotplus b)^k := (a+b)(a+qb)\cdots(a+q^{k-1}b)$
and consider the modified exponential function
$$h(z,a,b,q) = \sum_{k\geq 0} q^{\binom{k}{2}} \dfrac{(a \dotplus b)^k}{(1 \dotdiv q)^k}(-z)^k, \qquad (10)$$

which for $(a,b) = (0,1)$ reduces to $E_2(-z)$.
We want to study the series
$$f(z,a,b,q) = \dfrac{h(qz,a,b,q)}{h(z,a,b,q)}. \qquad (11)$$

From

$$h(z,a,b,q) - h(qz,a,b,q) = \sum_{k\geq 0} q^{\binom{k}{2}} \dfrac{(a \dotplus b)^k}{(1 \dotdiv q)^k}(-z)^k(1-q^k) =$$

$$= -z\sum_{k\geq 0} q^{\binom{k-1}{2}} \dfrac{(a \dotplus b)^{k-1}(a+q^{k-1}b)}{(1 \dotdiv q)^{k-1}}(-qz)^{k-1} = -azh(qz,a,b,q) - bzh(q^2 z,a,b,q)$$

we deduce that
$$f(z,a,b,q) = 1 + azf(z,a,b,q) + bzf(z,a,b,q)f(qz,a,b,q). \qquad (12)$$

If we write
$$f(z,a,b,q) = \sum_{n \geq 0} C_n(a,b,q) z^n, \tag{13}$$

then $C_n(a,b,q)$ is a $q-$analogue of $C_n(a,b) = \frac{1}{n} \sum_{k=1}^{n} \binom{n}{k}\binom{n}{k-1} b^{n-k}(a+b)^k$.

For it is easy to verify that the uniquely determined formal power series $f(z,a,b)$, which satisfies the functional equation $f(z,a,b) = 1 + azf(z,a,b) + bzf(z,a,b)^2$ has the series expansion

$$f(z,a,b) = 1 + \sum_{n \geq 1} \sum_{k=1}^{n} \frac{1}{n}\binom{n}{k}\binom{n}{k-1} b^{n-k}(a+b)^k z^n, \tag{14}$$

Of course $C_n(0,1) = \frac{1}{n+1}\binom{2n}{n} = C_n$ are the well known Catalan numbers.

The numbers $N(n,k) = \frac{1}{n}\binom{n}{k}\binom{n}{k-1}$ are called Narayana numbers. (cf. e.g. [10]).

Therefore we call $C_n(a,b,q)$ a $q-$Narayana polynomial. It satisfies

$$C_n(a,b,q) = aC_{n-1}(a,b,q) + b\sum_{k=0}^{n-1} q^k C_k(a,b,q) C_{n-k-1}(a,b,q) \tag{15}$$

with initial value $C_0(a,b,q) = 1$.

The first values are $1, a+b, a^2 + (2+q)ab + (1+q)b^2, \cdots$.

From the definition is clear that $C_n(a,b,q)$ is a polynomial in $a,b$ which is homogeneous of degree $n$. Therefore $C_n(a,b,q)$ has a unique representation in the form

$$C_n(a,b,q) = \sum_{k=0}^{n} N(n,k,q)(a+b)^k b^{n-k}. \tag{16}$$

There does not seem to be a simple formula for $N(n,k,q)$. All we can say is that they are polynomials in $q$ with integer coefficients and that $N(n,k,1) = N(n,k) = \frac{1}{n}\binom{n}{k}\binom{n}{k-1}$.

The first values of these $q-$Narayana numbers are given in the following table:

```
1
0    1
0    1        1
0    1        2 + q                1
0    1        3 + 3 q              3 + 2 q + q²           1
0    1        4 + 6 q + q² - q³    6 + 8 q + 5 q² + q³    4 + 3 q + 2 q² + q³    1
```

E.g. we get
$C_0(a,b,q) = 1,$
$C_1(a,b,q) = a+b,$
$C_2(a,b,q) = (a+b)b + (a+b)^2,$
$C_3(a,b,q) = (a+b)b^2 + (2+q)(a+b)^2 b + (a+b)^3.$

It is clear that $C_n(a,b,0) = (a+b)^n$. Therefore we have $N(n,k,0) = \binom{n-1}{k-1}$.

We shall need also another relation:

$$h(z,a,b,q) - h(qz,a,b,q) = \sum_{k \geq 0} q^{\binom{k}{2}} \frac{(a+b)^k}{(1 \div q)^k}(-z)^k(1-q^k) =$$

$$= -(a+b)z \sum_{k \geq 0} q^{\binom{k-1}{2}} \frac{(a+qb)^{k-1}}{(1 \div q)^{k-1}}(-qz)^{k-1} = -(a+b)zh(qz,a,qb,q).$$

This implies

$$\frac{f(qz,a,b,q) - 1}{f(z,a,b,q) - 1} = \frac{\frac{h(q^2z,a,b,q) - h(qz,a,b,q)}{h(qz,a,b,q)}}{\frac{h(qz,a,b,q) - h(z,a,b,q)}{h(z,a,b,q)}} = q \frac{f(qz,a,qb,q)}{f(z,a,b,q)}. \tag{17}$$

It is well known (cf. e.g. [10]) that the Catalan numbers $C_n = \frac{1}{n+1}\binom{2n}{n}$ can be characterized by the values of their Hankel determinants. They satisfy $\det(C_{i+j})_{i,j=0}^n = 1$ and $\det(C_{i+j+1})_{i,j=0}^n = 1$.

We want to show that the $q-$ Narayana polynomials can also be characterized by the values of their Hankel determinants. This shows that these polynomials are in some sense a natural generalization of the $q-$ Catalan numbers.

**Theorem**
*The polynomials $C_n(a,b,q)$ are characterized by their Hankel determinants*

$$\det(C_{i+j}(a,b,q))\big|_{i,j=0}^n = q^{\frac{n^2(n+1)}{2}} b^{\binom{n+1}{2}} (a+b)^n (a+qb)^{n-1} (a+q^2b)^{n-2} \cdots (a+q^{n-1}b) \tag{18}$$

*and*

$$\det(C_{i+j+1}(a,b,q))\big|_{i,j=0}^n = q^{\frac{n(n+1)^2}{2}} b^{\binom{n+1}{2}} (a+b)^{n+1} (a+qb)^n (a+q^2b)^{n-1} \cdots (a+q^nb). \tag{19}$$

**Remark**
*For the $q-$ Catalan numbers this reduces to*

$$\det(C_{i+j}(q))\big|_{i,j=0}^n = q^{\frac{n(n+1)(4n-1)}{6}} \quad \text{and} \quad \det(C_{i+j+1}(q))\big|_{i,j=0}^n = q^{\frac{n(n+1)(4n+5)}{6}}.$$

These special cases have been proved by another method in [5].

In order to prove the theorem we need the following well-known (cf. e.g. [1] or [8])

**Lemma**
*Let*

$$\sum_{k \geq 0} \mu_k z^k = \cfrac{1}{1 - s_0 z - \cfrac{t_0 z^2}{1 - s_1 z - \cfrac{t_1 z^2}{1 - \cdots}}} = \frac{1}{1 - s_0 z -} \frac{t_0 z^2}{1 - s_1 z -} \frac{t_1 z^2}{1 - s_2 z -} \cdots. \tag{20}$$

*Then the Hankel determinants have the following values*
$$\det(\mu_{i+j})_{i,j=0}^n = t_0^n t_1^{n-1} \cdots t_{n-1}$$
*and*
$$\det(\mu_{i+j+1})_{i,j=0}^n = d_{n+1} t_0^n t_1^{n-1} \cdots t_{n-1},$$
*where* $d_n = s_{n-1} d_{n-1} - t_{n-2} d_{n-2}, d_0 = 1, d_1 = s_0.$
*Let* $F$ *be the linear functional defined by* $F(z^n) = \mu_n.$
*Then the polynomials* $p_n(z),$ *defined by*
$p_0(z) = 1, p_1(z) = z - s_0,$
*and* $p_k(z) = (z - s_{k-1}) p_{k-1}(z) - t_{k-2} p_{k-2}(z)$
*satisfy* $F(p_n p_m) = t_0 \cdots t_{n-1} [n = m],$
*i.e. are orthogonal with respect to the linear functional* $F.$

**Proof of the Theorem**
We set
$$f(z, a, b, q) = 1 + (a+b) z g(\frac{z}{q}, a, b, q). \tag{21}$$

Then we can write (12)
in the form
$$f(z, a, b, q) = 1 + (a+b) z f(z, a, b, q) + q(a+b) b z^2 f(z, a, b, q) g(z, a, b, q), \tag{22}$$

where the series $g(z, a, b, q)$ satisfies
$$g(z, a, b, q) = 1 + q(a+b) z g(z, a, b, q) + q^2 b z g(qz, a, b, q) + q^3 b(a+b) z^2 g(z, a, b, q) g(qz, a, b, q).$$

We want to show that this is equivalent with

$$g(z, a, b, q) = 1 + q(a + b + qb) z g(z, a, b, q) + q^4 b(a + qb) z^2 g(z, a, b, q) g(qz, a, qb, q). \tag{23}$$

This means that
$q^2 b z g(qz, a, b, q) + q^3 b(a+b) z^2 g(z, a, b, q) g(qz, a, b, q) =$
$= q b z g(z, a, b, q) + q^4 b(a + qb) z^2 g(z, a, b, q) g(qz, a, qb, q)$
or
$g(qz, a, b, q)(1 + q(a+b) z g(z, a, b, q)) = g(z, a, b, q)(1 + (a + qb) q^2 z g(qz, a, qb, q)).$

By the definition of $g(z,a,b,q)$ this reduces to
$$g(qz,a,b,q)f(qz,a,b,q) = g(z,a,b,q)f(q^2z,a,qb,q),$$
which is equivalent with
$$\frac{g(qz,a,b,q)}{g(z,a,b,q)} = \frac{1}{q}\frac{(a+b)q^2zg(qz,a,b,q)}{(a+b)qzg(z,a,b,q)} = \frac{f(q^2z,a,qb,q)}{f(qz,a,b,q)}.$$
But this is an immediate consequence of (17).

Now (23) is equivalent with
$$g(z,a,b,q) = \frac{1}{1 - s_1 z - t_1 z^2 g(qz,a,qb,q)}$$
with $s_1 = q(a+b+qb)$ and $t_1 = q^4 b(a+qb)$.

This gives us a representation of $g(z,a,b,q)$ as a continued fraction of the form (20) with
$s_n = q^n(a + q^{n-1}b + q^n b)$ and $t_n = q^{3n+1}b(q^n b + a)$.
From (22) we conclude that $f(z,a,b,q)$ has also a representation as a continued fraction of the form (20) with the same $s_n, t_n$ together with $s_0 = a+b$ and $t_0 = q(a+b)b$.
From this the theorem immediately follows.
It is clear that the sequence $\big(C_n(a,b,q)\big)$ is uniquely determined by these determinants.

**Remark:** *It can be shown that the corresponding orthogonal polynomials are*
$$p_n(z,a,b,q) = \sum_{k=0}^{n}(-1)^{n-k}q^{\binom{n-k}{2}}z^k \sum_{j=k}^{n} q^{\binom{n+1}{2}-\binom{n+k+1-j}{2}} \begin{bmatrix} n+k-j \\ k \end{bmatrix}\begin{bmatrix} j-1 \\ k-1 \end{bmatrix} b^{j-k}(a+b)^{n-j}$$
*and*
$$p_n(z,0,1,q) = \sum_{k=0}^{n}(-1)^{n-k} q^{2\binom{n-k}{2}} \begin{bmatrix} n+k \\ 2k \end{bmatrix} z^k.$$
*But we don't need this result.*

c) Let now
$$h^*(z,a,b,q) := \sum_{k \geq 0} q^{\binom{k}{2}} \frac{r_k(a,b)}{(1-q)^k}(-z)^k, \qquad (24)$$

where $r_n(a,b,q) = \sum_{k=0}^{n}\begin{bmatrix} n \\ k \end{bmatrix} a^k b^{n-k}$ is a Rogers-Szegö polynomial, which satisfies
$r_n(a,b,q) = (a+b)r_{n-1}(a,b,q) + ab(q^{n-1}-1)r_{n-2}(a,b,q)$ (cf. e.g. [3]).

Therefore we get
$h^*(qz,a,b,q) - h^*(z,a,b,q) = -z(a+b)h^*(qz,a,b,q) - qabz^2 h^*(q^2z,a,b,q).$

If we define
$$f^*(z,a,b,q) = \frac{h^*(qz,a,b,q)}{h^*(z,a,b,q)}, \tag{25}$$

we see that $f^*(z,a,b,q)$ satisfies the functional equation
$$f^*(z,a,b,q) = 1 + (a+b)zf^*(z,a,b,q) + qabz^2 f^*(z,a,b,q)f^*(qz,a,b,q). \tag{26}$$

It is easy to see (compare (14) and (23)) that the series $f^*(z,a,b) = g(z, a-b, b, 1)$
satisfying $f^*(z,a,b) = 1 + (a+b)zf^*(z,a,b) + abz^2 f^*(z,a,b)^2$
has the expansion
$$azf^*(z,a,b) = \sum_{n\geq 1}\left(\sum_{k=1}^{n} N(n,k)a^k b^{n-k}\right)z^n.$$

Therefore we write $af^*(z,a,b,q) = \sum_{n\geq 0} C^*_{n+1}(a,b,q)z^n$

with
$$C^*_n(a,b,q) = \sum_{k=1}^{n} N^*(n,k,q)a^k b^{n-k}. \tag{27}$$

Again we have $N^*(n,k,0) = \binom{n-1}{k-1}$ and $D^*_n(1,1,0) = 2^{n-1}$.

The first values of $(N^*(n,k,q))_{k=1}^n, n \geq 1,$ are

```
1
1      1
1      2 + q              1
1      3 + 2 q + q²       3 + 2 q + q²              1
1      4 + 3 q + 2 q² + q³   6 + 6 q + 5 q² + 2 q³ + q⁴   4 + 3 q + 2 q² + q³   1
```

Since $f^*(z,a,b) = f^*(z,b,a)$ we see that $N^*(n,k,q) = N^*(n, n-k+1, q)$.
From this we get
$C^*_1(a,b,q) = a$
$C^*_2(a,b,q) = ab + a^2,$
$C^*_3(a,b,q) = ab^2 + (2+q)a^2 b + a^3,$
$C^*_4(a,b,q) = ab^3 + (3+2q+q^2)a^2 b^2 + (3+2q+q^2)a^3 b + a^4, \cdots.$

Let now $F(z,a,b,q) = 1 + azf^*(z,a,b,q)$. Then by (26)

$$F(z,a,b,q) = 1 + azF(z,a,b,q) - bzF(qz,a,b,q) + bzF(z,a,b,q)F(qz,a,b,q) \tag{28}$$

If we set $C^*_0(a,b,q) = 1$, then (28) implies the recurrence

$$C^*_n(a,b,q) = aC^*_{n-1}(a,b,q) + b\sum_{k=0}^{n-2} q^k C^*_k(a,b,q)C^*_{n-1-k}(a,b,q). \tag{29}$$

Comparing with [7] (5.5) we see that $C_n^*(1,s,q)$ are the Pólya-Gessel $q-$Catalan numbers, the first values of which are
$1,1,1+s,1+2s+qs+s^2,1+3s+2qs+q^2s+3s^2+2qs^2+q^2s^2+s^3,\cdots$

By (26) the corresponding $s_k,t_k$ are given by
$s_k=q^k(a+b)$ and $t_k=q^{2k+1}ab$. Therefore their Hankel determinants are

$$\det(C_{i+j+1}^*(a,b,q))_{i,j=0}^n=(ab)^{\binom{n+1}{2}}q^{\sum_{i=0}^n i^2}=(ab)^{\binom{n+1}{2}}q^{\frac{n(n+1)(2n+1)}{6}} \qquad (30)$$

and

$$\det(C_{i+j+2}^*(a,b,q))_{i,j=0}^n=(abq)^{\binom{n+1}{2}}q^{\frac{n(n+1)(2n+1)}{6}}\frac{a^{n+2}-b^{n+2}}{a-b}. \qquad (31)$$

It is easy to verify that the corresponding orthogonal polynomials are

$$p_n(z,a,b)=\sum_{k=0}^n(-1)^{n-k}q^{\binom{n-k}{2}}z^k\sum_{j=k}^n\begin{bmatrix}n+k-j\\k\end{bmatrix}\begin{bmatrix}j\\k\end{bmatrix}a^{j-k}b^{n-j}.$$

In order to compute

$\det(C_{i+j}^*(a,b,q))_{i,j=0}^n$

we observe that from (28) we get
$$F(z,a,b,q)=\frac{1-bzF(qz,a,b,q)}{1-bzF(qz,a,b,q)-az}=\frac{1}{1-\dfrac{az}{1-bzF(qz,a,b,q)}}=\frac{1}{1-}\frac{az}{1-}\frac{bz}{1-}\frac{qaz}{1-}\cdots=$$
$$=\frac{1}{1-azF(z,b,qa,q)}$$

Therefore we have
$F(z,a,b,q)=1+azF(z,a,b,q)F(z,b,qa,q)$
or
$F(z,a,b,q)=1+azF(z,a,b,q)+abz^2F(z,a,b,q)f^*(z,b,qa,q)$
This gives us $t_k=q^{2k}ab$ for all $k\geq 0$ and therefore we have

$$\det(C_{i+j}^*(a,b,q))_{i,j=0}^n=(ab)^{\binom{n+1}{2}}q^{\frac{n(n+1)(n-1)}{3}}.$$

In this case the orthogonal polynomials are

$$p_n(z,a,b)=\sum_{k=0}^n(-1)^{n-k}q^{\binom{n-k}{2}}z^k\sum_{j=k}^n\begin{bmatrix}n+k-j\\k\end{bmatrix}\begin{bmatrix}j-1\\j-k\end{bmatrix}b^{j-k}a^{n-j}.$$

For the special case $(a,b)=(1,s)$ these results have been proved by other methods in [6].
The generating function $f(z,q)$ of the $q-$Catalan numbers $C_n(q)$ has the continued fraction expansion
$$f(z,q)=\frac{1}{1-zf(qz,q)}=\frac{1}{1-\dfrac{z}{1-qzf(q^2z,q)}}=\frac{1-qzf(q^2z,q)}{1-qzf(q^2z,q)-z}.$$

Therefore we get
$$f(z,q) = 1 + zf(z,q) - qzf(q^2z,q) + qzf(z,q)f(q^2z,q).$$
This implies $f(z,q) = F(z,1,q,q^2)$ and thus the well known (cf. [7]) result that
$C_n(q) = C_n^*(1,q,q^2).$

This can also directly be seen:
$$h^*(z,1,q,q^2) = \sum_{k=0}^{n} q^{2\binom{k}{2}} \left( \sum_{j=0}^{k} \begin{bmatrix} k \\ j \end{bmatrix}_{q^2} q^k \right) \frac{(-z)^k}{(1 \div q^2)^k} = E_2(-z),$$
since
$$\sum_{j=0}^{k} \begin{bmatrix} k \\ j \end{bmatrix}_{q^2} q^k = (1 + q)^k \text{ (cf. [3])}.$$
From Gauss's formula $r_{2n+1}(1,-1) = 0, r_{2n}(1,-1) = (1-q)(1-q^3)\cdots(1-q^{2n-1})$ (cf e.g.[3])
we conclude in the same way that
$C_{2n+1}^*(1,-1,q) = (-1)^n q^n C_n(q^2)$
and $C_{2n+2}^*(1,-1,q) = 0.$

d) Another interesting special case is given by the $q-$Motzkin numbers $M_n(q)$ which have been considered in [5]. Their generating function $M(z) = \sum_{n \geq 0} M_n(q) z^n$ satisfies

$$M(z) = 1 + zM(z) + qz^2 M(z) M(qz). \text{ Therefore } M_n(q) = C_n^*(\frac{1+\sqrt{-3}}{2}, \frac{1-\sqrt{-3}}{2}, q).$$

The first values are
$1, 1, 1+q, 1+2q+q^2, 1+3q+3q^2+q^3+q^4, \cdots.$

The Hankel determinants are easily seen to be
$$\det(M_{i+j}(q))_{i,j=0}^{n} = q^{\frac{n(n+1)(2n+1)}{6}}$$
and
$$\det(M_{i+j+1}(q))_{i,j=0}^{n} = q^{\frac{n(n+1)(2n+1)}{6} + \binom{n+1}{2}} d_{n+1} = q^{2\binom{n+2}{3}} d_{n+1}.$$
Here $(d_n)_{n \geq 0} = (1,1,0,-1,-1,0,1,1,0,\cdots)$ is periodic with period 6.